\input amssym.def
\input amssym
\magnification=1200
\parindent0pt
\hsize=16 true cm
\baselineskip=13  pt plus .2pt
$ $

\def\G{(\bar \Gamma,{\cal G})}

\def \S {{\cal S}}
\def\B{{\Bbb B}}
\def\S{{\Bbb S}}

\centerline {\bf On large groups of symmetries of finite graphs embedded in spheres}

\bigskip \bigskip

\centerline {Bruno P. Zimmermann}

\medskip

\centerline {Universit\`a degli Studi di Trieste}

\centerline {Dipartimento di Matematica e Geoscienze}

\centerline {34127 Trieste, Italy}

\bigskip  \bigskip

{\bf Abstract.}  Let $G$ be a finite group acting orthogonally on a pair
$(S^d,\Gamma)$ where $\Gamma$ is a finite, connected graph of genus $g>1$
embedded in the sphere $S^d$. The 3-dimensional case $d=3$ has
recently been considered in a paper by C. Wang, S. Wang, Y. Zhang and the present
author where for each genus $g>1$ the maximum order of a $G$-action on a pair
$(S^3,\Gamma)$ is determined and the corresponding graphs $\Gamma$ are classified. 
In the present paper we consider arbitrary dimensions $d$ and prove  that the order of
$G$ is bounded above by a polynomial of degree $d/2$ in $g$ if $d$ is even and of
degree $(d+1)/2$ if $d$ is odd; moreover the degree $d/2$ is best possible in even
dimensions $d$. We discuss also the problem, given a finite graph $\Gamma$ and its finite
symmetry group, to find the minimal dimension of a sphere into which $\Gamma$ embeds
equivariantly as above.

\bigskip


\medskip


\bigskip \bigskip

{\bf 1. Introduction}

\medskip

We study large finite groups $G$ of automorphisms of a finite,
connected  graph $\Gamma$ which embeds smoothly into a sphere $S^d$ of some dimension
$d$ such that the $G$-action on $\Gamma$ extends to an orthogonal action of
$G$ on $S^d$. In other words,  we study large finite groups
$G$ of orthogonal transformations of  pairs  $(S^d, \Gamma)$ where $\Gamma$ denotes a
finite, connected graph smoothly embedded in a sphere $S^d$. All actions considered in
the present paper will be faithful on both $S^d$ and $\Gamma$, and all finite graphs
$\Gamma$ will be {\it hyperbolic}, i.e. connected,  of genus $g>1$ (the rank of its
free fundamental group) and {\it without free edges} (edges with one vertex of valence
1; note that free edges can be deleted in a
$G$-equivariant way without changing the genus of a graph); we allow closed and
multiple edges.

\medskip

The case of dimension $d = 3$ is considered in [6]. For a finite subgroup of
the orthogonal group SO(4) acting on a pair 
$(S^3, \Gamma)$, a regular neighbourhood of $\Gamma$ in $S^3$ is a 3-dimensional
handlebody $V^3_g$ of genus $g>1$ on which $G$ acts orientation-preservingly, and by
[8], [1] there is the linear bound $|G| \le 12(g-1)$ for orientation-preserving actions
on 3-dimensional handlebodies.  In [6], for each genus
$g>1$ the maximal possible order of a $G$-action on $(S^3, \Gamma)$ is determined and the
corresponding graphs $\Gamma$ are classified (on the basis of analogous results in [5]
for the case of closed surfaces embedded in $S^3$); the maximal possible order
$12(g-1)$ is obtained only for finitely many values of $g$.

\medskip

Concerning dimension $d=4$, supppose that $G$ acts orthogonally on a pair $(S^4,
\Gamma)$; now a regular neighbourhood of $\Gamma$ is a 4-dimensional handlebody $V^4_g$
whose boundary $\partial V^3_g$ is a connected sum $\sharp_g (S^2 \times S^1)$ of $g$
copies of $S^2 \times S^1$. Finite group actions on such connected sums are considered
in [11] whose results imply the quadratic upper bound $|G| \le 24g(g-1)$, for $g
\ge 15$; moreover there does not exist a linear bound in $g$ for the order of $G$.

\medskip

In the following main result of the present paper  we consider the case of
arbitrary dimensions $d \ge 4$.

\bigskip

{\bf Theorem 1.}  {\sl  Let $G$ be a finite subgroup of the orthogonal group ${\rm
O}(d+1)$ acting on a pair $(S^d, \Gamma)$, for a finite hyperbolic graph of genus
$g > 1$ embedded in $S^d$.  Then the order of $G$ is bounded above by a polynomial
of degree $d/2$ in $g$ if $d$ is even and of degree
$(d+1)/2$ if $d$ is odd. The degree $d/2$ is best possible in even dimensions  whereas
in odd dimensions the optimal degree is either $(d-1)/2$ or $(d+1)/2$.}

\bigskip

So the optimal degree in odd dimensions remains open at present (except for $d=3$ where
it is $(d-1)/2 = 1$). The proof of Theorem 1 will be reduced to an analogous result in
[2] about finite group actions on $d$-dimensional handlebodies  (see Theorem 2 in the
next section).

\medskip

The maximum order of a finite group of automorphisms of a finite
hyperbolic graph of genus $g > 2$ is $2^gg!$ ([7]), obtained for a graph with one vertex
and $g$ closed edges (a "bouquet of $g$ circles") whose automorphism group is isomorphic
to the semidirect product $(\Bbb Z_2)^g \ltimes S_g$ (with normal subgroup 
$(\Bbb Z_2)^g$ on which the symmetric group $S_g$ acts by permutation of
coordinates).  At present we don't know the minimal dimension of a sphere which admits an
equivariant embedding of the bouquet of $g$ circles (i.e., invariant under an
orthogonal action of $(\Bbb Z_2)^g \ltimes S_g$), see the question at the end of section
2. For some other graphs with large symmetry groups instead, we determine this minimal
dimension of an equivariant embedding in the examples in section 2.

\bigskip

{\bf 2. Proof of Theorem 1}

\medskip

See [2] for the following.  A $d$-dimensional handlebody $V^d_g$ of genus $g$ can be
defined as a regular neighbourhood of a finite connected graph $\Gamma$ of genus g
embedded in $S^d$. Such a handlebody $V^d_g$ can be uniformized by a Schottky group
${\cal S}_g$, a free group of rank $g$ of M\"obius transformations of $S^{d-1}$ which
extends naturally to the disk $\B^d$ ("Poincar\'e extension"); the interior of
$B^d$ is the Poincar\'e-model of hyperbolic space $\Bbb H^d$ on which 
M\"obius transformations act as hyperbolic isometries. The handlebody $V^d_g$ is
obtained as the quotient $(B^d - \Lambda({\cal S}_g))/{\cal S}_g$ where
$\Lambda({\cal S}_g)  \subset S^{d-1}$ denotes the set of limit point of the action of
${\cal S}_g$ on $B^d$ (a Cantor set), in particular
$B^d - \Lambda({\cal S}_g)$ is the universal covering of 
$V^d_g$.  This gives the interior $\Bbb H^d/{\cal S}_g$ of $V^d_g$ the structure of a
complete hyperbolic manifold, and we say that the Schottky group ${\cal S}_g$
uniformizes the {\it hyperbolic handlebody} $V^d_g$. In particular there is the notion of
an isometry of such a hyperbolic handlebody meaning that it acts as an isometry on
the interior of $V^d_g$; equivalently, each lift to the universal covering $B^d -
\Lambda({\cal S}_g)$ of $V^d_g$ extends to a M\"obius transformation of $B^d$.

\medskip

The following is proved in [2].

\bigskip

{\bf Theorem 2.}  {\sl  Let $G$ be a finite group of isometries of a hyperbolic
handlebody $V^d_g$ of dimension $d \ge 3$ and of genus $g > 1$ which acts faithfully
on the fundamental group. Then the order of $G$ is bounded by a polynomial of degree
$d/2$ in $g$ if $d$ is even and of degree
$(d+1)/2$ if $d$ is odd.}

\bigskip

Since a handlebody of dimension $d \ge 4$ admits  $S^1$-actions there is
no upper bound for the order of finite group actions which are not
faithful on the fundamental group (however there is a Jordan-type bound for such
actions, see [2, Corollary]). On the other hand, finite faithful actions on finite
hyperbolic graphs are faithful on the fundamental group, that is inject
into the outer automorphism group of the fundamental group  ([12, Lemma 1]); 
conversely, it is observed in [9, p. 478] (as a version of the Nielsen realization
problem for free groups) that  every finite subgroup $G$ of the outer automorphism group
of a free group can be realized by an action of $G$ on a finite graph.

\bigskip

Starting with the {\it Proof of Theorem 1} now,  let $G$ be a finite group acting
orthogonally on $(S^d, \Gamma)$, for a finite hyperbolic graph $\Gamma$ of genus
$g>1$ embedded in $S^d$. A $G$-invariant regular neighbourhood of $\Gamma$ in
$S^d$ is homeomorpic to a handlebody $V^d_g$ of dimension $d$ and genus $g$. Since we
are assuming that the action of $G$ on the graph $\Gamma$ is faithful, by [12, Lemma 1]
also the induced action of $G$ on the fundamental group of $\Gamma$ and hence of
$V^d_g$ is faithful and defines an injection of $G$ into the outer automorphism
groups of the fundamental groups of $\Gamma$ and
$V^d_g$. The first part of Theorem 1 is now a consequence of Theorem 2 and the
following:

\bigskip

{\bf Proposition.}  {\sl  Let $G$ be a finite group acting orthogonally on a pair
$(S^d, \Gamma)$, for a finite hyperbolic graph of genus
$g>1$ embedded in $S^d$, and hence also on a handlebody $V^d_g$ in $S^d$ obtained as a
$G$-invariant regular neighbourhood of $\Gamma$. Then $V^d_g$ can be uniformized by a
Schottky group such that $V^d_g$ admits an isometric
$G$-action (inducing the same injection into the outer automorphism group of the
fundamental group as the original $G$-action).}

\bigskip

Very likely the original $G$-action on
$V^d_g$ is in fact conjugate to an isometric $G$-action; however in order to
apply Theorem 2 we just need some isometric $G$-action on $V^d_g$ which is faithful on
the fundamental group, so we don't follow this here.

\bigskip

{\it Proof of the Proposition.}  The group $G$ acts as a group of automorphisms of the
finite graph $\Gamma$.  Let $\tilde \Gamma$ be the universal covering tree of $\Gamma$;
the group of all liftings of all elements of $G$ to $\tilde \Gamma$ defines a group $E$ of
automorphisms of the tree $\tilde \Gamma$ and a group extension 
$1 \to F_g  \hookrightarrow  E  \to G \to 1$ where $F_g$ denotes the universal covering
group, a  free group of rank $g$ isomophic to the fundamental group of $\Gamma$.

\medskip

By possibly subdividing edges by a new vertex  
we can assume that $G$ acts on $\Gamma$ without inversions of edges; then the quotient
$\bar \Gamma = \Gamma/G$ is again a finite graph. Choose a maximal tree in $\bar \Gamma$
and lift it isomorphically first to
$\Gamma$ and then to $\tilde \Gamma$, then lift also the remaining edges of $\bar
\Gamma$ to $\Gamma$ and $\tilde \Gamma$. Associating to the vertices and edges of 
$\bar \Gamma$ the stabilizers in $G$ or $E$ of the lifted vertices and edges in 
$\Gamma$ and $\tilde \Gamma$, this defines a finite graph of finite groups $\G$, with
inclusions of the edge groups into the adjacent vertex groups. The fundamental group
$\pi_1\G$ of the finite graph of finite groups $\G$ is the iterated free product with
amalgamation and HNN-extension of the vertex groups amalgamated over the edge groups,
first taking the iterated free product with amalgamation over the chosen maximal tree of
$\Gamma$ and then associating an HNN-generator to each of the remaining edges. By the
standard theory of groups acting on trees, graphs of groups and their fundamental
groups (see [3], [4] or [10]), the extension $E$ is isomorphic to  
the fundamental group $\pi_1\G$ of $\G$ and we have a group extension
$$1 \to F_g  \hookrightarrow  E = \pi_1\G  \to G \to 1.$$

We will assume in the following that the graph of groups $\G$ has no {\it trivial
edges}, i.e. no edges with two different vertices such that the edge group coincides
with one of the two vertex groups (by collapsing trivial edges, i.e. by amalgamating the
two vertex groups into a single vertex group). We will realize the extension 
$E = \pi_1\G$ as a group of isometries of hyperbolic space $\Bbb H^d$, the interior of the
$d$-ball $B^d$, or equivalently as a group of M\"obius transformations of $B^d$ or
$S^{d-1} = \partial B^d$ such that the subgroup $F_g$ is realized by a Schottky group
${\cal S}_g$. Then  $G$ acts as a group of isometries of the hyperbolic handlebody 
$V^d_g = (B^d - \Lambda({\cal S}_g))/{\cal S}_g$ proving the Proposition.
The realization of $E = \pi_1\G$ is by standard combination methods as 
described in [2, p.247] and [9, p. 479-482]. As an illustration we discuss the
case of a graph of groups $\G$ with a single, non-closed edge. We lift the edge to an
edge $B$ of $\Gamma$, with vertices
$A_1$ and $A_2$. The corresponding edge and vertex groups of $\G$ (which we denote by
the same letters) are defined as the stabilizers of the edge $B$ and its vertices in the
group $G$ acting orthogonally on $\Gamma \subset S^d$, and  
$E = \pi_1\G$ is the free product with amalgamation $A_1 *_BA_2$.

\medskip

For $i = 1$ and 2, the vertex group $A_i \subset G$ acts orthogonally on $S^d$ and fixes
the corresponding vertex $A_i$ of the edge $B$ of $\Gamma$. A  regular
invariant neighbourhood of the fixed point $A_i$ is a ball
$B_i^d$ with an action of $A_i$. The edge group $B$ is a subgroup of
both  $A_1$ and $A_2$ which fixes both vertices $A_1$ and $A_2$ and pointwise the
edge $B$,. We can assume that the intersection $\partial B_1^d \cap \partial B_2^d$ 
of the two boundary spheres $S^{d-1}$ is non-empty and hence a sphere $S^{d-2}$. We
conjugate $A_2$ by the reflection in this sphere $S^{d-2}$ and obtain a group
$A_2'$ fixing the vertex $A_1$; since the reflection commutes with each element of 
$B$, both $A_1$ and $A'_2$ act on $B_1^d$ now with a common subgroup $B$.

\medskip

Identifying $B_1^d$ with the standard ball $B^d$, with the center $A_1$
corresponding to the origin 0, we have orthgonal actions of $A_1$ and $A_2'$ on $B^d$
with a common subgroup $B$. The orthogonal action of $B$ on $B^d$ fixes two diametral
points in $\partial B^d$ which are not fixed by any other element of $A_1$ and $A'_2$
(corresponding to the intersection of the edge $B$ with $\partial B_1^d$), and hence $B$
fixes pointwise the hyperbolic line $L$ in hyperbolic space $\Bbb H^d$ (the interior of
$B^d$) connecting these two diametral points.  We conjugate $A'_2$ by a reflection in a
hyperbolic hyperplane orthogonal of the line $L$ far from the origin $0 \in B^d$. The
group generated by $A_1$ and the reflected group
$A''_2$ is then isomorphic to the free product with amalgamation $A_1 *_BA''_2$ and
realizes 
$E = A_1 *_BA_2$ as a group of hyperbolic isometries of $\Bbb H^d$, or 
equivalently of M\"obius transformations of $S^{d-1} = \partial B^d$ (by standard
combination methods, see [2, p. 247] and [9, p. 480]
for more details and some figures).

\medskip

In a similar way, inductively edge by edge in finitely many steps,  one realizes the
group $E = \pi_1\G$ also in the general case (see again [9, p. 479-482]).

\bigskip

This concludes the proof of the Proposition
and, by Theorem 2, also of the first statement of Theorem 1.

\bigskip

For the second statement of Theorem 1, we construct an infinite series of orthgonal
actions of finite groups $G$ on finite graphs $\Gamma$ embedded in $S^d$ which realize
the lower bounds for the polynomial degrees in Theorem 1.

\bigskip

{\bf Example 1.}  For $k > 1$, let $G = C_1 \times \ldots \times C_k  \cong (\Bbb
Z_m)^k$, of order $n = m^k$, be the product of $k$ cyclic groups $C_i \cong \Bbb Z_m$
of order
$m$. Choose an orthogonal action of $G$ on  $\Bbb R^{2k}$ as follows. Decomposing
$\Bbb R^{2k} = P_1 \times \ldots \times P_k$ as the product of $k$ orthogonal planes
$P_i$, each $C_i$ acts on $P_i$ faithfully by rotations and trivially on the $k-1$
orthogonal planes. This $G$-action on $\Bbb R^{2k}$ extends to an orthogonal $G$-action
on the one-point compactification $S^{2k}$ of $\Bbb R^{2k}$, with two global fixed
points 0 and $\infty$.

\medskip

We consider a graph $\Gamma$ in $S^{2k}$ with two vertices 0 and $\infty$ and $km$
connecting edges divided into $k$ groups of $m$ edges. We embed the first
group of $m$ edges into $P_1$ such that $C_1$ permutes these edges cyclically,
then the next $m$ edges into $P_2$ etc., defining an orthogonal action of
$G$ on the graph $\Gamma$ embedded in $S^{2k}$. The graph $\Gamma$ has genus $g = km - 1$,
hence 
$$|G| = m^k = (g + 1)^k/k^k$$ which is a polynomial of degree $k=d/2$ in $g$.

\medskip

Suppose that $m > 2$; then $d = 2k$ is the minimal dimension of a sphere which admits a
$G$-equivariant embedding of the graph $\Gamma$. In fact, in such an embedding into a
sphere $S^d$ the group  $G \cong (\Bbb Z_m)^k$ has a global fixed point and hence acts
orthogonally on the boundary $S^{d-1}$ of an invariant regular neighbourhood of a
global fixed point, and $d-1 = 2k-1$ is the minimal dimension of a sphere with a
faithful orthogonal action of $(\Bbb Z_m)^k$; equivalently, the minimal dimension of a
faithful, real, linear representation of $(\Bbb Z_m)^k$ is $2k$.

\medskip

In odd dimensions $d = 2k+1$, we extend the orthogonal action of $G$ on
$\Bbb R^{2k}$ to an orthogonal action on
$\Bbb R^{2k+1}$, trivial on the last coordinate, and then proceed as
before; we get a polynomial of degree $k = (d-1)/2$ in
$g$ for the order of $G$. As noted in the Introduction, the optimal degree in
dimension $d = 3$ is $(d-1)/2 = 1$; for odd dimensions  $d > 3$ the optimal
degree is either $(d-1)/2$ or $(d+1)/2$ but at present we do not know which of these
two values occurs.

\medskip

This completes the proof of Theorem 1.

\bigskip

We present some other infinite series  of finite orthogonal group actions on finite
graphs embedded in spheres which realize the minimal dimension of such an embedding.

\bigskip

{\bf  Examples 2.}  i)  Let $\Gamma$ be the complete graph with $d+2$ vertices, or
the the 1-skeleton of a $(d+1)$-simplex. We embed the regular $(d+1)$-simplex into
$B^{d+1}$, with vertices in $S^d = \partial B^{d+1}$, and project its edges
radially to  $S^d$; this defines an embedding of $\Gamma$ into $S^d$. The
automorphism group of the graph $\Gamma$ is the symmetric group $S_{d+2}$ which extends to
an orthogonal action on the sphere $S^d$. Again $d$ is  the minimal dimension of such an
equivariant embedding since it is the minimal dimension of a sphere with a
faithful, orthogonal action of  $S_{d+2}$ (equivalently, the minimal dimension of a
faithful, real, linear representation of the symmetric group $S_{d+2}$ is $d+1$).

\medskip

ii)  Let $\Gamma'$ be a graph with two vertices
and $d+2$ egdes. Then $\Gamma'$ has an embedding into $S^{d+1}$ 
which is "dual" to the embedding of
$\Gamma$ into $S^d \subset S^{d+1}$ in part i) in the following sence : the sphere $S^d$
seperates $S^{d+1}$ into two balls $B^{d+1}$; the two vertices of $\Gamma'$ are the
centers of these two balls, and each of the $d+2$ connecting edges intersects exactly one
of the $d+2$ faces of the projected $(d+1)$-simplex in $S^d$ in its center.

\medskip

The symmetric group $S_{d+2}$ acts orthogonally on
$(S^{d+1}, \Gamma')$ with two global fixed points (the double suspension of the action of
$S_{d+2}$ on $S^d$), in particular it acts on
$S^d$ and it follows as in i) that $d+1$ is the minimal dimension of a sphere with an
equivariant embedding of $\Gamma'$.

\medskip

iii)  Let $\Gamma$ be the 1-skeleton of the 
$(d+1)$-dimensional hypercube now, with a projection to $S^d = \partial B^{d+1}$ as in i).
The graph $\Gamma$ has an automorphism group $(\Bbb Z_2)^{d+1}$  which extends to an 
orthogonal action on $S^d$ (inversion of coordinates), and $d$ is the minimal
dimension of a sphere with an orthogonal action of  $(\Bbb Z_2)^{d+1}$ (the minimal
dimension of a faithful, real, linear representation of $(\Bbb Z_2)^{d+1}$ is $d+1$).

\medskip

Dualizing as in ii) we get a graph $\Gamma'$ in $S^{d+1}$ with two vertices connected by
$2(d+1)$ edges (each intersecting exactly one of the faces of the hypercube), with an
action of  $(\Bbb Z_2)^{d+1}$ with two global fixed
points, and
$d+1$ is the minimal dimension of such an embedding (since $d$ is the minimal dimension
of an orthgonal action of $(\Bbb Z_2)^{d+1}$ on a sphere). Note that this
realizes the minimal dimension in the case $m=2$ of the group $({\Bbb Z_2})^k$ in
Example 1, replacing each plane $P_i$ with a rotation by a line with a reflection.

\bigskip

{\bf Question.}  As noted in the introduction, the maximum order of a finite group of
automorphisms of a finite hyperbolic graph of genus $g>2$ is $2^gg!$ ([7]),
obtained for a graph $\Gamma$ with a single vertex and $g$ closed edges whose 
automorphism group is the semidirect product $(\Bbb Z_2)^g \ltimes S_g$ (the symmetric
group $S_g$ acts by permutation of the edges, the nomal subgroup  $(\Bbb Z_2)^g$ 
by inversion of the edges).  What is the minimal dimension of a sphere with an equivariant
embedding of this graph?

\smallskip

Note that the graph $\Gamma$ embeds equivariantly into $\Bbb R^2 \times \ldots \times \Bbb
R^2 = \Bbb R^{2g}$, by embedding each of the $g$ closed edges into a different plane
$\Bbb R^2$ (the unique vertex corresponds to the origin of $\Bbb R^{2g}$ and also of each
plane
$\Bbb R^2$); the group $(\Bbb Z_2)^g$ acts by reflections in lines through the origin on
the $g$ planes, inverting the $g$ embedded closed edges, the symmetric group $S_g$ acts by
permutation of the planes.  Then
$\Gamma$ admits an equivariant embedding also into  the one-point
compactification $\S^{2g}$ of $\Bbb R^{2g}$. On the
other hand, the minimal dimension of a sphere with an orthgonal action of $(\Bbb Z_2)^g
\ltimes S_g$ is
$g-1$ (since $g$ is the minimal dimension of a faithful, real, linear representation of
$(\Bbb Z_2)^g$).  It seems reasonable, however, that the minimal dimension of an
equivariant embedding of $\Gamma$ into a sphere is $2g$ (as described above), but this
minimal dimension, between $g-1$ and $2g$, remains open at present.

\bigskip \bigskip

\centerline {\bf References}

\bigskip

\item {[1]} D. McCullough, A. Miller, B. Zimmermann,  Group actions on
handlebodies,  {\it Proc. London Math. Soc.}  59  (1989)  373-415

\item {[2]}  M. Mecchia, B. Zimmermann,  On finite groups of isometries of
handlebodies in arbitrary dimensions and finite extensions of Schottky groups, 
{\it  Fund. Math.} 230 (2015)  237-249

\item {[3]} P. Scott, T. Wall,  Topological methods in group theory, 
{\it Homological Group Theory,  London Math. Soc. Lecture Notes} 36, {\it
Cambridge University Press}  1979

\item {[4]} J.P. Serre,  Trees,   {\it Springer} 1980

\item {[5]}  C. Wang, S. Wang, Y. Zhang, B. Zimmermann,  Embedding
surfaces into $S^3$ with maximum symmetry,  
{\it Groups Geometry and Dynamics} 9  (2015)  1001-1045

\item {[6]}  C. Wang, S. Wang, Y. Zhang, B. Zimmermann, Graphs in the
3-sphere with maximum symmetry,  arXiv:1510.00822

\item {[7]} S. Wang, B. Zimmermann, The maximum order finite groups of outer
automorphisms of free groups,  {\it Math. Z.}  216  (1994)  83-87

\item {[8]} B. Zimmermann,  \"Uber Abbildungsklassen von Henkelk\"orpern,  {\it Arch.
Math.} 33  (1979)  379-382

\item {[9]} B. Zimmermann,  \"Uber Hom\"oomorphismen n-dimensionaler
Henkelk\"orper und endliche Erweiterungen von Schottky-Gruppen,  {\it Comm. Math. Helv.}
56   (1981)  474-486

\item {[10]} B. Zimmermann,  Generators and relations for discontinuous groups, 
{\it  Generators and Relations in Groups and Geometries (eds. Barlotti, Ellers, Plaumann,
Strambach),  NATO Advanced Study Institute Series} 333,  {\it Kluwer Academic
Publishers}  (1991)  407-436

\item {[11]}   B. Zimmermann,  On finite groups acting on a connected sum of
3-manifolds $S^2 \times S^1$, {\it  Fund. Math.} 226 (2014) 131-142

\item {[12]} B. Zimmermann,   Finite groups of outer automorphism groups of free
groups  {\it  Glasgow Math. J.} 38  (1996)  275-282

\bye